\documentclass[]{article}
\usepackage{graphicx,subfigure}
\usepackage{amsfonts,amsbsy,amscd,amsgen,amsthm,dsfont,amsmath,amssymb,mathrsfs}
\usepackage[colorlinks=true,allcolors=blue]{hyperref}
\usepackage{bm,multirow,esvect}
\theoremstyle{definition}
\usepackage{tabularx}
\usepackage{enumitem}
\usepackage[utf8]{inputenc}
\usepackage[T1]{fontenc}
\usepackage{authblk}
\usepackage{algorithm,algpseudocode}
\usepackage{textcomp}
\usepackage{xcolor}
\usepackage{hhline}

\newcommand{\bxi}{\pmb \xi}
\newcommand{\vc}{\mathbf}

\renewcommand{\tilde}{\widetilde}

\newcommand{\Null}{\mathcal{N}}
\newcommand{\A}{\mathcal{A}}

\DeclareMathAlphabet\mathbfcal{OMS}{cmsy}{b}{n}
\DeclareMathOperator*{\argmin}{arg\,min}

\newcommand{\err}{\mathcal E}
\newcommand{\bphi}{\pmb\phi}
\newcommand{\balpha}{\pmb\alpha}

\newcommand{\R}{\mathcal{R}}

\renewcommand{\span}{\mbox{span}}

\begin{document}

\title{State Estimation Using Sparse DEIM and Recurrent Neural Networks} 
\author[]{Mohammad Farazmand\thanks{Corresponding author's email address: farazmand@ncsu.edu}}
\affil{Department of Mathematics, North Carolina State University, 2311 Stinson Drive, Raleigh, NC 27695-8205, USA}
\date{}
\maketitle

\begin{abstract}
Sparse Discrete Empirical Interpolation Method (S-DEIM)	was recently proposed for state estimation in dynamical systems when  only a sparse subset of the state variables can be observed. The S-DEIM estimate involves a kernel vector whose optimal value is inferred through a data assimilation algorithm. This data assimilation step suffers from two drawbacks: (i) It requires the knowledge of the governing equations of the dynamical system, and (ii) It is not generally guaranteed to converge to the optimal kernel vector. To address these issues, here we introduce an equation-free S-DEIM framework that estimates the optimal kernel vector from sparse observational time series using recurrent neural networks (RNNs). We show that the recurrent architecture is necessary since the kernel vector cannot be estimated from instantaneous observations. But RNNs, which incorporate the past history of the observations in the learning process, lead to nearly optimal estimations. We demonstrate the efficacy of our method on three numerical examples with increasing degree of spatiotemporal complexity: a conceptual model of atmospheric flow known as the Lorenz-96 system, the Kuramoto--Sivashinsky equation, and the Rayleigh--B\'enard convection. In each case, the resulting S-DEIM estimates are satisfactory even when a relatively simple RNN architecture, namely the reservoir computing network, is used. More specifically, our RNN-based S-DEIM state estimations reduce the relative error between 42\% and 58\% when compared to Q-DEIM which ignores the kernel vector by setting it equal to zero.
\end{abstract}

\section{Introduction}\label{sec:intro}
Estimating the state of a dynamical system from incomplete observations is essential since it enables downstream tasks such as forecasting and control~\cite{law2015,Thau1973}. As such, state estimation has widespread application in robotics, fluid dynamics, meteorology, oceanography, and many other fields~\cite{Barfoot2017,Callaham2019}. 

Among existing state estimation techniques, Discrete Empirical Interpolation Method (DEIM) has quickly gained traction over the past decade for its simplicity and interpretability. DEIM was originally developed for reduced-order modeling of dynamical systems~\cite{Sorensen2010}.
It was later recognized that DEIM can also be used to estimate the state of a system from its incomplete pointwise observations~\cite{Manohar2018}. These observations are gathered by sensors. The optimal locations of the sensors are often approximated by the column-pivoted QR decomposition, resulting in the so-called Q-DEIM method~\cite{Drmac2016}.

As the number of observations increase, the (Q-)DEIM state estimation improves. However, when relatively few sensors are available, the estimation error can be very large. Farazmand~\cite{Farazmand2024} recently proposed Sparse DEIM (S-DEIM) to reduce the estimation error, even when few observations are available. To this end, S-DEIM considers the general solution to the least squares problem at the center of DEIM, as opposed to the minimum-norm solution used by previous DEIM-based methods.
This general solution leads to the estimation, 
\begin{equation}
 \vc u(t) \simeq \underbrace{\Phi \Theta^+\vc y(t)}_{\text{closed form}} + \underbrace{\Phi \vc z(t)}_{\text{learned}},
 \label{eq:SDEIM_gen}
\end{equation}
where $\vc u(t)$ is the state to be estimated, $\vc y(t)$ represents the observations at time $t$, $\Phi$ and $\Theta$ are appropriate matrices as defined in Section~\ref{sec:prelim}, and the superscript $+$ denotes the pseudo inverse of a matrix. The first term in~\eqref{eq:SDEIM_gen} can immediately be computed from the observations $\vc y(t)$. The second term involves the kernel vector $\vc z(t)$ which belongs to the null space of $\Theta$ and whose optimal value needs to be estimated. The (Q)-DEIM estimate uses to the trivial kernel vector $\vc z(t) = \vc 0$. In contrast, the S-DEIM estimate uses a nonzero optimal kernel vector.

Unfortunately, the optimal kernel vector cannot be inferred from instantaneous observational data $\vc y(t)$. Instead, Farazmand~\cite{Farazmand2024} introduced a data assimilation method which seeks to approximate this vector using the governing equations of the dynamical system. There are two main issues with this data assimilation approach: (i) Although the method converges to the optimal kernel vector under certain assumptions, it is not guaranteed to converge generally. (ii) It requires the full knowledge and use of the differential equations which govern the dynamical systems. Thus, if the governing equations are unknown, the data assimilation becomes infeasible. Even when the governing equations are known, the data assimilation step can be computationally expensive.

To address these issues, here we develop an equation-free method to approximate the optimal kernel vector in S-DEIM. Our method relies on a Recurrent Neural Network (RNN) to learn the optimal kernel vector $\vc z(t)$ from time series of the observational data $\vc y(s)$ with $s\leq t$, leveraging the information encoded in the past history of the observations. As such, the main innovation in this paper is the use of RNNs to estimate the optimal kernel vector in S-DEIM, without relying on the governing equations of the dynamical systems. We demonstrate, with several numerical examples, that the resulting state estimations are nearly optimal.

\subsection{Related work}
In this section, we review the historical development of DEIM-based methods and their applications.

\textit{Model order reduction:} 
Projection-based reduced models seek to approximate a high-dimensional dynamical system with a smaller system which evolves on a linear subspace of the state space~\cite{Holmes1996,holmes97}. If this linear subspace approximates the attractor of the original system, one hopes that the reduced model estimates the dynamics at a lower computational cost. In this context, the reduced equations are traditionally derived by orthogonal projection of the dynamics on the reduced subspace~\cite{Benner2015}. This approach has a computational complexity of $\mathcal O(mN)$ where $N$ is the dimension of the original system and $m\ll N$ is the dimension of the reduced subspace~\cite{Greif2019,Babaee2023}. 
As a result, the computational cost of evaluating the reduced model can be similar to that of the original system, negating the purpose of model order reduction.

To address this issue, Chaturantabut and Sorensen~\cite{Sorensen2010} proposed DEIM which replaces the orthogonal projection to the reduced subspace with an oblique projection; also see~\cite{Barrault2004}. The computational complexity of DEIM is $\mathcal O(m)$, which significantly reduces the computational cost when compared to orthogonal projection methods.

Several improvements and generalizations of DEIM have been proposed. For instance, Drma{\v c} and Saibaba~\cite{Drmac2018} generalize DEIM to weighted inner product spaces. DEIM can be sensitive to noise which prompted its stability analysis in~\cite{Argaud2017,Drmac2020}. In \cite{Saibaba2020}, randomized range finding and randomized subset selection algorithms were developed to further reduce the computational cost of DEIM. Esmaeili et al.~\cite{Esmaeili2020} developed Generalized DEIM (G-DEIM) which identifies an optimal reduced subspace by leveraging the Jacobian information in addition to the state variables. Peherstorfer et al.~\cite{Willcox2014} proposed Localized DEIM (L-DEIM) which, instead of a global subspace, uses a collection of subspaces, each locally approximating the attractor. Farazmand and Saibaba~\cite{Farazmand2023} developed a tensor-based variant of DEIM which significantly reduces the storage cost.

\textit{State estimation:} 
As the name suggests, DEIM can be used for interpolation purposes. In the context of dynamical systems, this interpolation amounts to state estimation from incomplete observations~\cite{Mohseni2014,Manohar2018,Callaham2019}. Let the state variable be $N$-dimensional and assume that only an $r$-dimensional subset of them can be observed. As we review in Section~\ref{sec:prelim}, DEIM can be used to approximate the remaining $N-r$ variables within a reduced subspace spanned by $m$ modes. 

In fact, the overdetermined case where $r>m$ predates DEIM. It was first developed by Everson and Sirovich~\cite{Sirovich1995} in the context of image reconstruction. Later, it was dubbed as \emph{gappy proper orthogonal decomposition} (or gappy POD) and was applied to dynamical systems~\cite{Willcox2004,Willcox2006}.

The underdetermined case, where $r<m$, has received relatively little attention, in spite of its relevance in applications where gathering observational data is burdensome~\cite{Brunton2021}. Numerical experiments suggest that the DEIM estimation error can be large in the underdetermined case~\cite{Manohar2018,kakasenko2025}. Motivated by these observations, Farazmand~\cite{Farazmand2024} developed Sparse DEIM (S-DEIM) to reduce the estimation error. To this end, S-DEIM takes advantage of a kernel vector which had been ignored in the previous formulations of DEIM. The present paper is concerned with this underdetermined regime ($r<m$). We use recurrent neural networks, which leverage the temporal history of sparse observations, to estimate the optimal kernel vector in S-DEIM.

\textit{Sensor placement:} 
The quality of DEIM-based state estimation depends significantly on the location of the sensors which gather the observational data. For instance, if the sensor locations are assigned at random, the DEIM estimation error is typically very large. Drma{\v c} and Gugercin~\cite{Drmac2016} proposed the column-pivoted QR (CPQR) factorization as a greedy algorithm for approximating the set of optimal sensor locations. The resulting method is dubbed Q-DEIM and is reviewed in Section~\ref{sec:CPQR} below. In the overdetermined case where $r>m$, CPQR only returns the optimal locations for the first $m$ sensors. Hendryx et al.~\cite{Hendryx2021} proposed extended DEIM (E-DEIM) to approximate the optimal location of the remaining $r-m$ sensors; also see~\cite{Shin2024}.

Although CPQR is most commonly used for sensor placement in DEIM, alternative methods can also be employed. For instance, Saito et al.~\cite{Saito2021} proposed a determinant-based method which in some cases outperforms the CPQR algorithm.
Klishin et al.~\cite{Klishin2023} use a similar determinant-based method which also takes into account the pairwise interactions between sensors. 

\subsection{Outline}
The remainder of this paper is organized as follows. In Section~\ref{sec:prelim}, we review the mathematical preliminaries concerning S-DEIM. In Section~\ref{sec:RNN}, we introduce a machine learning framework which estimates the unknown optimal kernel vector in S-DEIM using recurrent neural networks. Section~\ref{sec:num_results} contains our numerical results. Finally, we present our concluding remarks in Section~\ref{sec:conc}

\section{Sparse DEIM}\label{sec:prelim}
We consider a dynamical system whose state at time $t$ is denoted by $\vc u(t)\in\mathbb R^N$. The state can, for instance, denote the solution of an ordinary differential equation (ODE) or a spatially discretized solution of a time-dependent partial differential equation (PDE). Table~\ref{tab:glossary} summarizes the key terms and notation used in this paper.

\renewcommand{\arraystretch}{1.2}
\begin{table}
	\centering
	\caption{Glossary of key terms and notation.}
	\begin{tabular}{|c |l |}
		\hline
		\textbf{Symbol} & \textbf{Description} \\ \hhline{|=|=|}
		$N$ & State Space Dimension \\ \hline
		$r$ & Number of Sensors \\ \hline
		$m$ & Number of Modes (Basis Functions) \\ \hline
		$\vc u(t)\in\mathbb R^N$ & Unknown State\\ \hline
		$\vc y(t)\in\mathbb R^r$ & Known Observations\\ \hline
		$S_r \in\mathbb R^{N\times r}$ & Selection Matrix, $\vc y(t)=S_r^\top\vc u(t)$\\ \hline
		$\Phi_m = [\pmb \phi_1|\cdots|\pmb \phi_m]\in\mathbb R^{N\times m}$ & Basis Matrix\\ \hline
		$\tilde{\vc u}(\vc z(t))\in\mathbb R^N$ & S-DEIM Estimate\\ \hline
		$\hat{\vc u}(t)\in\mathbb R^N$ & Best Fit, $\hat{\vc u}(t)=\Phi_m\Phi_m^\top \vc u(t)$\\ \hline
		$\hat{\vc z}(t)\in\mathbb R^m$ & Optimal Kernel Vector,  $\hat{\vc z}(t)=ZZ^\top \Phi_m^\top \vc u(t)$\\ \hline
		$\mathcal E_m(\balpha):=\|\Phi_m \balpha - \vc u\|$ & Truncation Error\\ \hline
		$\err_{obs}(\balpha)  = \|S_r^\top\Phi_m\balpha -\vc y\|$ & Observation Residual\\ \hline
	\end{tabular}
	\label{tab:glossary}
\end{table}

Assuming that only a sparse subset of the entries of $\vc u(t)$ are known through observations, we would like to estimate the missing entries. More precisely, let $\{i_1,i_2,\cdots,i_r\}\subset \{1,2,\cdots,N\}$ denote the distinct indices of the observed state variables $\{u_{i_1}, u_{i_2},\cdots,u_{i_r}\}$. We denote the observations by a vector $\vc y = [u_{i_1},\cdots,u_{i_r}]^\top+\pmb\epsilon$ where $\pmb\epsilon\in\mathbb R^r$ is the observational noise.
Defining the \emph{selection matrix} $S_r=\left[\vc e_{i_1}|\vc e_{i_2}|\cdots|\vc e_{i_r}\right]\in\mathbb R^{N\times r}$, the observations can be written as
\begin{equation}
\vc y(t) = S_r^\top \vc u(t) + \pmb\epsilon(t).
\end{equation}
Given the observations $\vc y(t)$, we would like to estimate the entire state $\vc u(t)$. For the theoretical review in this section, we ignore the observational noise and set $\pmb\epsilon=\vc 0$. 
In Section~\ref{sec:num_results}, where we present our numerical results, we add a $5\%$ Gaussian noise to the observations to demonstrate that the results are robust to noise.

Consider an approximation of the state $\vc u(t)$ within an orthonormal set $\left\{\bphi_1,\cdots,\bphi_m\right\}\subset \mathbb R^N$ so that 
\begin{equation}
	\vc u(t) \simeq \overline{\vc u}+ \sum_{i=1}^m \alpha_i(t) \bphi_i=\overline{\vc u}+ \Phi_m\balpha(t),
\end{equation}
where $\overline{\vc u}$ is the mean, $\Phi_m = [\bphi_1|\bphi_2|\cdots|\bphi_m]\in\mathbb R^{N\times m}$ is the \emph{basis matrix}, and $\balpha = [\alpha_1,\alpha_2,\cdots,\alpha_m]^\top\in\mathbb R^m$ are the coordinates in that basis. Without loss of generality, we assume $\overline{\vc u}=\vc 0$. The \emph{best fit approximation} is given by
\begin{equation}
	\hat{\vc u} := \Phi_m \Phi_m^\top \vc u,
	\label{eq:bestFit}
\end{equation}
corresponding to the coefficients $\hat \balpha = \Phi_m^\top \vc u$ which minimize the truncation error $\mathcal E_m(\balpha):=\|\Phi_m \balpha - \vc u\|$, where $\|\cdot\|$ denotes the standard Euclidean norm of a vector. The truncation error $\mathcal E_m$ is a monotonically decreasing function of $m$, i.e., using more modes leads to better approximations.
Unfortunately, however, the best fit approximation~\eqref{eq:bestFit} is not computable since it relies on the full state $\vc u$.

To overcome this issue, DEIM determines the coefficients $\balpha$ by minimizing the discrepancy between the observations $\vc y$ and the reconstructed state at the sensor locations. To this end, we define the observation residual, 
\begin{equation}
\err_{obs}(\balpha)  = \|S_r^\top\Phi_m\balpha -\vc y\|.
\label{eq:err_obs}
\end{equation}
Note that $S_r^\top\Phi_m\balpha$ corresponds to the entries of the reconstruction $\Phi_m\balpha$ at the sensor locations.
To obtain the optimal coefficients, we minimize $\err_{obs}$ over all $\balpha$. The general solution to
this optimization problem is given by~\cite{Farazmand2024}
\begin{equation}
(S_r^\top \Phi_m)^+\vc y+\vc z=\argmin_{\balpha\in\mathbb R^m} \err_{obs}(\balpha) , 
\end{equation}
for an arbitrary $\vc z\in\Null[S_r^\top \Phi_m]$,
 where the superscript $+$ denotes the Moore-Penrose pseudo-inverse and $\Null$ denotes the null space of a linear operator. The matrix $\Theta$ in equation~\eqref{eq:SDEIM_gen} is equal to $S_r^\top \Phi_m$.

\begin{figure}
\centering
\includegraphics[width=.6\textwidth]{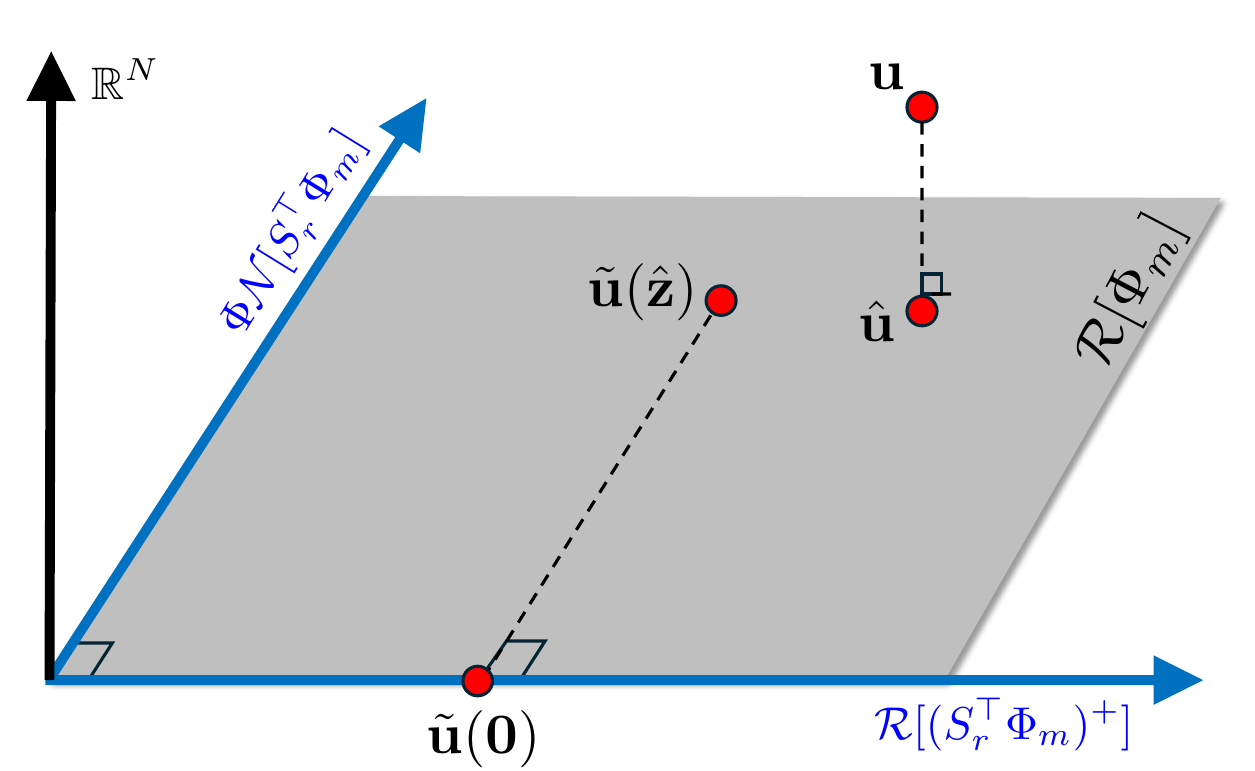}
\caption{State space geometry of DEIM estimation with $r<m$. The truncated basis $\{\bphi_1,\cdots,\bphi_m\}$ spans the subspace $\R[\Phi_m]\subset\mathbb R^N$. The variables denote the true state $\vc u$, the best fit $\hat{\vc u}$, DEIM estimate $\tilde{\vc u}(\vc 0)$, 
and the optimal S-DEIM estimate $\tilde{\vc u}(\hat{\vc z})$. The subspace $\R[\Phi_m]$ admits the orthogonal decomposition $\R[\Phi_m]=\R[(S_r^\top\Phi_m)^+]\oplus \Phi\Null[S_r^\top\Phi_m]$.}
\label{fig:SDEIM_geom}
\end{figure}

The corresponding S-DEIM estimate is given by
\begin{equation}
\tilde{\vc u}(\vc z) = \Phi_m(S_r^\top \Phi_m)^+\vc y+\Phi_m \vc z,\quad \vc z\in \Null[S_r^\top\Phi_m].
\label{eq:SDEIM}
\end{equation}
(Q)-DEIM uses the trivial kernel vector $\vc z=\vc 0$, leading to the state estimation $\tilde{\vc u}(\vc  0)$. Figure~\ref{fig:SDEIM_geom} shows a schematic depiction of the state space geometry of the true state, the best fit, DEIM and S-DEIM estimations.
The DEIM estimate $\tilde{\vc u}(\vc 0)$ is computable because it only requires the observations $\vc y$ and not the entire state $\vc u$. In fact, in the overdetermined case, where the number of sensors is at least equal to the number of modes ($r\geq m$), we have $\Null[S_r^\top\Phi_m]=\{\vc 0\}$. In this case, the DEIM estimation is the only choice. 

However, in the underdetermined case, where the number of available sensors is limited so that $r<m$, the null space $\Null[S_r^\top\Phi_m]$ is non-trivial and therefore nonzero choices of the kernel vector $\vc z$ are possible. This underdetermined regime is the focus of the present paper. It arises naturally in problems where the number of available sensors $r$ is small, but $r$ basis functions are not sufficient to accurately estimate the state of the system, hence demanding the use of a larger set of basis vectors.

In this underdetermined case, (Q)-DEIM estimate, which uses the trivial kernel vector $\vc z=\vc 0$, is quite restrictive. As illustrated in figure~\ref{fig:SDEIM_geom}, the DEIM estimate $\tilde{\vc u}(\vc 0)$ belongs to the range of $(S_r^\top \Phi_m)^+$, which is a lower-dimensional subspace of $\mathcal R[\Phi_m]$. As such, the DEIM estimate is constrained to $\mathcal R[(S_r^\top \Phi_m)^+]$ and cannot take advantage of the full expressivity of the basis $\Phi_m$. On the other hand, the S-DEIM estimates, which use a non-zero $\vc z$, are free of this constraint and can explore the full range of the basis matrix $\Phi_m$ more freely.

The natural question then is whether there exists an optimal kernel vector $\vc z\in\Null[S_r^\top\Phi_m]$ such that the corresponding reconstruction $\tilde{\vc u}(\vc z)$ is closest to the true state $\vc u$. To answer this question, we minimize the total error, 
\begin{equation}
	\err_{tot}(\vc z) = \|\tilde{\vc u}(\vc z) - \vc u\|.
\end{equation}
We define the \emph{optimal kernel vector} $\hat{\vc z}$ as the kernel vector that minimizes the total error, 
\begin{equation}
	\hat{\vc z} = \argmin_{\vc z\in\Null[S_r^\top\Phi_m]}\err_{tot}(\vc z).
\end{equation}
Let $\{\vc z_1,\cdots,\vc z_{m-r}\}$ be an orthonormal basis for $\Null[S_r^\top\Phi_m]$ and define the \emph{kernel matrix} $Z=[\vc z_1|\cdots|\vc z_{m-r}]\in\mathbb R^{m\times (m-r)}$.
It is straightforward to show that the optimal kernel vector is given by~\cite{Farazmand2024}
\begin{equation}\label{eq:opt_z}
	\hat{\vc z}=ZZ^\top \Phi_m^\top \vc u.
\end{equation}
However, this optimal kernel vector is not directly computable from observations; it requires the knowledge of the full state $\vc u$.
In Section~\ref{sec:RNN}, we introduce a machine learning method for estimating the optimal kernel vector from observational time series.
But first we briefly review the S-DEIM error estimate (Section~\ref{sec:SDEIM_err}) as well as the common computational methods for determining the basis matrix $\Phi_m$ (Section~\ref{sec:POD}) and the selection matrix $S_r$ (Section~\ref{sec:CPQR}).

\subsection{S-DEIM error}\label{sec:SDEIM_err}
Let $\tilde{\vc u}(\vc z)$ denote the S-DEIM estimate~\eqref{eq:SDEIM} with an arbitrary kernel vector $\vc z\in \Null[S_r^\top\Phi_m]$. We would like to quantify the estimation error $\| \tilde{\vc u}(\vc z)-\vc u\|$. It was shown in~\cite[Theorem 1]{Farazmand2024} that this error can be decomposed into three terms, 
\begin{equation}\label{eq:SDEIM_err}
	\| \tilde{\vc u}(\vc z)-\vc u\|^2 = \| \vc u-\hat{\vc u}\|^2 + \| \mathbb D(\vc u-\hat{\vc u})\|^2+\|\vc z-\hat{\vc z}\|^2,\quad \forall \vc z\in \Null[S_r^\top\Phi_m],
\end{equation}
where $\mathbb D = \Phi_m (S_r^\top\Phi_m)^+S_r^\top$, $\hat{\vc u}=\Phi_m\Phi_m^\top\vc u$ is the best fit, and $\hat{\vc z}=ZZ^\top \Phi_m^\top \vc u$ is the optimal kernel vector~\eqref{eq:opt_z}.
We refer to~\cite{Farazmand2024} for the derivation of this identity.

The first term on the right-hand side of~\eqref{eq:SDEIM_err} involves the truncation error $\| \vc u-\hat{\vc u}\|$, reflecting the fact that the basis matrix $\Phi_m$ contains a truncated basis of $\mathbb R^N$. This error can be decreased by increasing the number of modes $m$.

The second term $\| \mathbb D(\vc u-\hat{\vc u})\|$ quantifies the effects of incomplete observational data. Given a fixed and limited number of sensors $r<m$, this component of the error can be decreased by identifying an optimal sensor placement strategy; see Section~\ref{sec:CPQR}.

Finally, the most relevant to our discussion is the third term $\|\vc z-\hat{\vc z}\|$ involving the kernel vector $\vc z$. Our main goal in this paper is to estimate and use the optimal kernel vector $\hat{\vc z}$ in order to minimize this component of the error as much as possible. (Q)-DEIM, which uses the trivial kernel $\vc z=\vc 0$, is bound to at least incur an error of $\|\hat{\vc z}\|$, no matter how many modes $m$ are used or how cleverly the sensors are placed.

\subsection{Basis matrix: Proper orthogonal decomposition}\label{sec:POD}
In principle, the basis matrix $\Phi_m = [\bphi_1|\bphi_2|\cdots|\bphi_m]$ can comprise any set of orthonormal vectors which is deemed suitable for approximating the state of the dynamical system. For instance, if the system is dissipative, the physically relevant dynamics take place on an attractor $\A\subset \mathbb R^N$.
In that case, the basis matrix $\Phi_m$ should be chosen such that its range $\R[\Phi_m]=\span\{\bphi_1,\cdots,\bphi_m\}$ approximates the attractor $\A$.

Such an approximating subspace is often obtained by the Proper Orthogonal Decomposition (POD)~\cite{Lumey1975,Sirovich1987}, also known as Principle Component Analysis (PCA)~\cite{Pearson1901,Hotelling1933}. 
POD requires samples from the attractor of the system. Let $\{\vc u_1, \vc u_2, \cdots, \vc u_M\}\subset \mathcal A$ denote $M$ such samples. In practice, these samples are obtained from carrying out high-fidelity numerical simulations or experiments. Defining the data matrix $A = [\vc u_1| \vc u_2|\cdots|\vc u_M]\in\mathbb R^{N\times M}$, we compute its Singular Value Decomposition (SVD),
\begin{equation}
A = \Phi \Sigma V^\top,
\end{equation}
where $\Phi = [\bphi_1|\cdots|\bphi_N]\in \mathbb R^{N\times N}$ is the left singular matrix, $\Sigma\in\mathbb R^{N\times M}$ contains the singular values of $A$, and $V\in\mathbb R^{M\times M}$ is the right singular matrix. The columns of the left singular matrix $\Phi$ are referred to as the \emph{POD modes}. Truncating these modes to the first $m$ columns, we obtain the basis matrix $\Phi_m$. It is known that $\R[\Phi_m]$ is the $m$-dimensional linear subspace that best approximates the snapshots $\{\vc u_1, \vc u_2, \cdots, \vc u_M\}$~\cite{Koch2007}.

In the above description of POD, we assumed that the empirical mean of the samples $\overline{\vc u} := \frac{1}{M}\sum_{i=1}^M \vc u_i$ vanishes. If the mean is nonzero, the same procedure can be applied to the centered data matrix $A = [\vc u_1-\overline{\vc u}|\cdots|\vc u_M-\overline{\vc u}]$ whose columns are the mean-removed samples~\cite{Rathinam2003}. Therefore, we can safely assume that the data has zero mean, i.e., $\overline{\vc u}=\vc 0$.

\subsection{Sensor placement: Column-pivoted QR decomposition}\label{sec:CPQR}
Drma\v{c} and Gugercin~\cite{Drmac2016} used the column-pivoted QR (CPQR) decomposition for sensor placement in DEIM. The resulting method is referred to as Q-DEIM. CPQR was first proposed for the special case $r=m$ in~\cite{Drmac2016} and later extended to the case with $r<m$ in~\cite{kakasenko2025}.
In the present work, we also use the CPQR algorithm for sensor placement, and therefore review the method briefly in this section.

To motivate CPQR, we first observe that the DEIM reconstruction $\tilde{\vc u}(\vc 0)$ satisfies following upper bound for its relative error~\cite{Farazmand2024}, 
\begin{equation}
	\frac{\|\tilde{\vc u}(\vc 0)-\vc u\|}{\|\vc u\|} \leq \|(S_r^\top \Phi_m)^+\|_2,
	\label{eq:ub}
\end{equation}
where $\|\cdot\|_2$ denotes the spectral norm (i.e., 2-norm) of a matrix. In the context of DEIM, sensor placement methods seek to determine the selection matrix $S_r$ that minimizes the upper bound $\|(S_r^\top \Phi_m)^+\|_2$.

Next, we recall that 
\begin{equation}
	\|(S_r^\top \Phi_m)^+\|_2 = \sigma_{max}((S_r^\top \Phi_m)^+) = \left[\sigma_{min}(S_r^\top \Phi_m) \right]^{-1},
\end{equation}
where $\sigma_{max}$ denotes the largest singular value of a matrix and $\sigma_{min}$ denotes its smallest nonzero singular value.
Therefore, minimizing the spectral norm $\|(S_r^\top \Phi_m)^+\|_2$ is equivalent to maximizing $\sigma_{min}(S_r^\top \Phi_m)$.

On the other hand, consider a rank-revealing QR decomposition of $\Phi_m^\top\in\mathbb R^{m\times N}$ with column pivoting, 
\begin{equation}
	\Phi_m^\top\Pi = Q R,\quad \Pi= [\Pi_1\quad \Pi_2],\quad R= \begin{bmatrix}
		R_{11} & R_{12}\\
		0 & R_{22}
	\end{bmatrix} 
	\label{eq:cpqr}
\end{equation}
where $\Pi$ is a permutation matrix with $\Pi_1\in\mathbb R^{N\times r}$ and $\Pi_2\in\mathbb R^{N\times (N-r)}$, $Q\in\mathbb R^{m\times m}$ is an orthogonal matrix, and $R_{11}\in\mathbb R^{r\times r}$ is an upper triangular matrix with diagonal entries $R_{11}^{ii}$. The CPQR algorithm takes $\Phi_m^\top$ as input and returns the matrices $(\Pi,Q,R)$ as outputs.
The permutation matrix $\Pi$ is determined such that $|R_{11}^{11}|\geq |R_{11}^{22}|\geq\cdots\geq |R_{11}^{rr}|$. Consequently, CPQR ensures that $\sigma_{min}(R_{11})$ is relatively large.
As shown in~\cite{kakasenko2025}, we have
\begin{equation}
\sigma_{min}(\Phi_m^\top \Pi_1) = \sigma_{min}(\Pi_1^\top \Phi_m) = \sigma_{min}(R_{11}).
\label{eq:sig_min}
\end{equation}
As a result, if we set $S_r=\Pi_1$, the upper bound~\eqref{eq:ub} is relatively small since $\|(S_r^\top \Phi_m)^+\|_2 = [\sigma_{min}(S_r^\top \Phi_m)]^{-1}=[\sigma_{min}(R_{11})]^{-1}$. In the following, we compute the CPQR decomposition of $\Phi_m^\top$ and set the selection matrix $S_r$ equal to the first $r$ columns of the permutation matrix $\Pi$.

When the selection matrix $S_r$ is determined by CPQR, the resulting matrix $S_r^\top\Phi_m$ has full row rank~\cite{gu1996}. Recall that the pseudo-inverse of a full row rank matrix is a right inverse, so that $(S_r^\top \Phi_m)(S_r^\top \Phi_m)^+ = I_r$. More specifically, if $A$ has full row rank, then $A^+= A^\top (AA^\top)^{-1}$.

We emphasize that CPQR sensor placement is derived based on the assumption that $\vc z=\vc 0$; see equation~\eqref{eq:ub}. As such, it is not guaranteed to return optimal (or even nearly optimal) sensor locations for S-DEIM state estimation which uses $\vc z\neq \vc 0$. Nonetheless, in the absence of a tailor-made sensor placement method for S-DEIM, we use CPQR here.

\section{S-DEIM Using Recurrent Neural Networks}\label{sec:RNN}
When $S_r^\top\Phi_m$ has full row rank, S-DEIM estimation~\eqref{eq:SDEIM} is an interpolation in the following sense. For any kernel vector $\vc z$, we have the \emph{interpolation property}~\cite{Farazmand2024},
\begin{equation}\label{eq:interp}
	S_r^\top \tilde{\vc u}(\vc z(t)) = \vc y(t),\quad \forall \vc z(t)\in\Null[S_r^\top\Phi_m].
\end{equation}
This implies that, regardless of the choice of the kernel vector, the estimation $\tilde{\vc u}(\vc z)$ agrees with the true observations $\vc y$.
As a result, the \emph{instantaneous} observations $\vc y(t)$ do not provide any additional information which could be used to estimate the optimal kernel vector $\hat{\vc z}(t)$~\cite{Farazmand2024}. 
\begin{figure}[h]
	\centering
	\includegraphics[width=\textwidth]{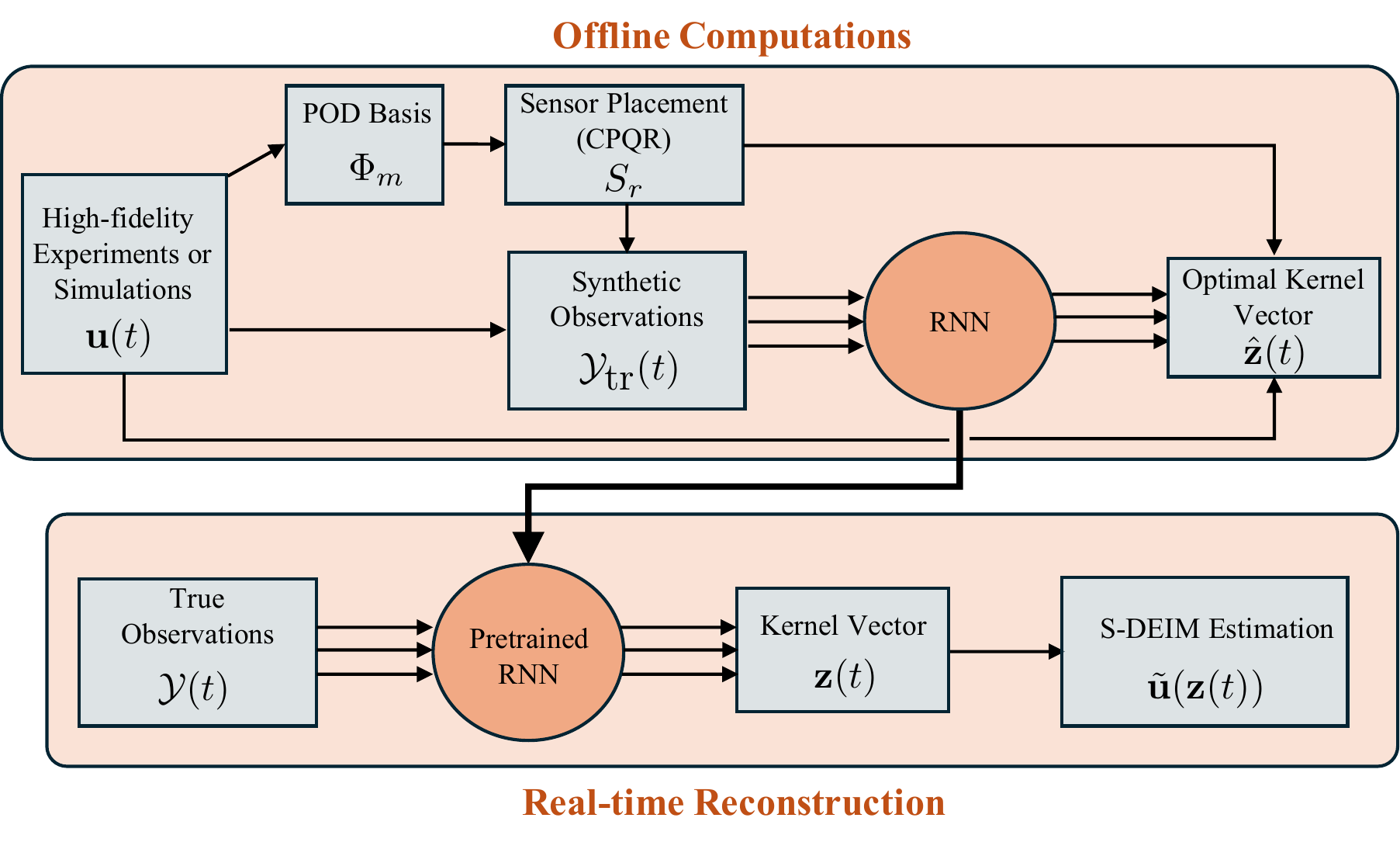}
	\caption{Pipeline of the proposed framework. In the offline stage, the POD modes and optimal sensor locations are computed and the RNN is trained. 
	In the real-time reconstruction, the pretrained network is used to estimate the state from observational data. The observational time series from the training set is denoted by $\mathcal Y_{tr}(t)$.}
	\label{fig:schem_pipeline}
\end{figure}

Instead, here we use the time series of the observations, including the past history, to estimate the optimal kernel vector at the current time $t$.
To this end, we define the time series,
\begin{equation}
	\mathcal Y(t) = \{\vc y(s): s\leq t\},
	\label{eq:ts}
\end{equation}
which contains all the observations in the past up until the current time $t$. As depicted in figure~\ref{fig:schem_pipeline}, these time series are used as input to a recurrent neural network in order to predict its output, i.e., the optimal kernel vector $\hat{\vc z}(t)$. The feasibility of such predictions is motivated by Takens' delay embedding theorem which states that the attractor of a dynamical system can be faithfully reconstructed from the time series of an observable~\cite{takens,sauer91}.

To train the RNN, we need the true outputs $\hat{\vc z}(t)$. The training is carries out in the offline computations part of the process using high-fidelity experiments or numerical simulations; see the top compartment of figure~\ref{fig:schem_pipeline}. Recall that such high-fidelity data is already required by (Q-)DEIM in order to compute the POD basis matrix $\Phi_m$ and 
to perform the CPQR decomposition to obtain the selection matrix $S_r$. Here, we leverage this available data to prepare the input $\mathcal Y(t)$ and the output $\hat{\vc z}(t)$ of the RNN.
Once the training is completed in the offline phase, the pretrained RNN is used to estimate the optimal kernel vector from actual sparse observational data without requiring any further high-fidelity data.  To this end, any RNN architecture can be used in principle. Here, we use the simplest form of RNNs, i.e., reservoir computing (RC) network~\cite{Jaeger2001,Maass2002}. 

The RC architecture is chosen here solely for its simplicity and low training cost. More complex networks, such as Long Short-Term Memory (LSTM) networks, can certainly be used instead and will likely improve the resulting S-DEIM estimates. Nonetheless, as we show with our numerical results (Section~\ref{sec:num_results}), RC-based S-DEIM already improves the state estimation accuracy by 42-58\%.

\subsection{Training the reservoir computing network}\label{sec:RC_training}
Before describing the RC architecture, we first outline the data required for its training. We assume access to a set of high-fidelity data $\vc u(t_i)$ at discrete time instances $t_0<t_1<t_2<\cdots<t_M$.
We further assume that the data is equispaced in time so that $t_{i+1}-t_i = \Delta t$. From this high-fidelity data set we compute the corresponding observations $\vc y(t_i) = S_r^\top \vc u(t_i)$ and the optimal kernel vectors $\hat{\vc z}(t_i)=ZZ^\top \Phi_m^\top \vc u(t_i)$.
Our aim is to train a neural network whose inputs are the observational time series $\vc y(t_i)$ and its outputs are the optimal kernel vectors $\hat{\vc z}(t_i)$. However, it is more prudent to train the neural networks with outputs being $\hat\bxi(t_i)=Z^\top \hat{\vc z}(t_i)=Z^\top \Phi_m^\top \vc u(t_i)$. Recall that the columns of the kernel matrix $Z$ form an orthonormal basis for $\Null[S_r^\top \Phi_m]$, and therefore the vector $\hat\bxi(t_i)\in\mathbb R^{m-r}$ denotes the coordinates of $\hat{\vc z}(t_i)\in\mathbb R^{m}$ in that basis. Since there is a one-to-one correspondence between $\hat{\vc z}(t_i)$ and $\hat\bxi(t_i)$, replacing the network output with $\hat\bxi(t_i)$ does not lead to any loss of information; the optimal kernel vector can be computed uniquely as $\hat{\vc z}(t_i)=Z\hat\bxi(t_i)$. Yet, the vector $\hat\bxi(t_i)$ has dimension $m-r$, which is smaller than the dimension of the optimal kernel vector $\hat{\vc z}(t_i)\in\mathbb R^m$. This lower dimensionality leads to a faster and more accurate training process. 
\begin{figure}
	\centering
	\includegraphics[width=0.9\textwidth]{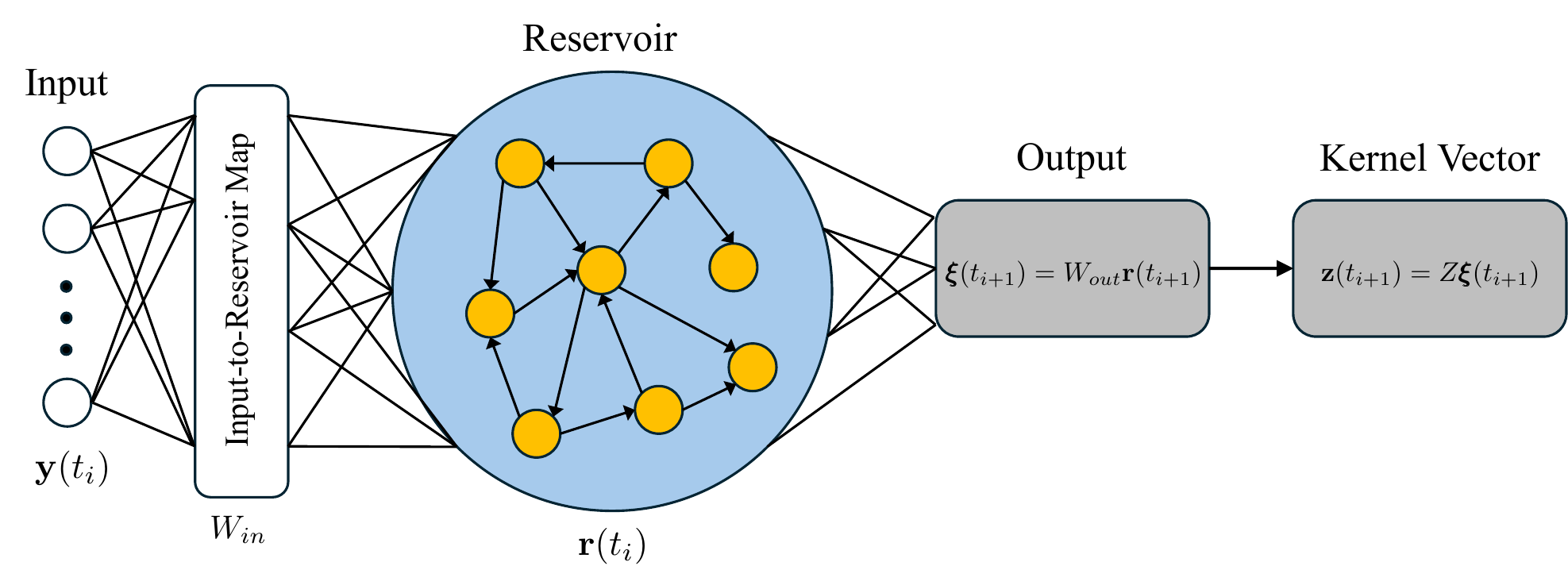}
	\caption{Architecture of the reservoir computing network.}
	\label{fig:RC_archit}
\end{figure}

Figure~\ref{fig:RC_archit} shows the architecture of our RC network. The observations $\vc y(t_i)$ constitute the inputs to the network. During the training phase, the outputs are the optimal kernel vectors $\hat{\bxi}(t_i)=Z^\top \Phi_m^\top \vc u(t_i)$, determined from high-fidelity states $\vc u(t_i)$. The reservoir state is denoted by $\vc r(t_i)\in\mathbb R^K$ where $K\gg \max\{r,m\}$. The inputs $\vc y(t_i)$ are connected to the reservoir state through the input-to-reservoir matrix $W_{in}\in\mathbb R^{K\times r}$. The nodes of the reservoir are connected to each other through a sparse reservoir weight matrix $W_R\in\mathbb R^{K\times K}$ and the biases $\vc b\in\mathbb R^K$. Finally, the reservoir state $\vc r(t_i)$ is mapped to the outputs through the output weight matrix $W_{out}\in\mathbb R^{(m-r)\times K}$.

What sets RC networks apart from most other neural networks is the fact that its input weights $W_{in}$, its reservoir weights $W_R$ and the biases $\vc b$ of the network are prescribed at random, as opposed to being determined through an optimization process. The only trained part of the network is the output weights $W_{out}$. As we will see shortly, there exists a simple closed form expression for the optimal output weights. As such, training RC networks is significantly faster than other recurrent neural networks such as LSTM networks.

At each time step, the reservoir state is updated according to the iterative map,
\begin{equation}\label{eq:res_state}
	\vc r(t_{i+1}) = (1-\alpha) \vc r(t_i) + \alpha \sigma(W_R\vc r(t_i) +W_{in}\vc y(t_i) +\vc b), 
\end{equation}
where $0<\alpha\leq 1$ is the learning rate and the nonlinear activation function $\sigma:\mathbb R\to\mathbb R$ is applied component-wise when its argument is a vector. Throughout this paper, we use a hyperbolic tangent activation funciton, so that $\sigma (x) = \tanh(x)$.

Ideally, the output weights $W_{out}$ are determined so that $\hat\bxi(t_i) = W_{out}\vc r(t_i)$. Since this might be an overdetermined problem, one instead solves the regularized least squares problem, 
\begin{equation}\label{eq:Wout_min}
	W_{out} = \argmin_{W} \|WR-\Xi\|_F^2+\lambda \|W\|_F^2,
\end{equation}
where $0<\lambda\ll 1$ is a regularization parameter, $\|\cdot\|_F$ denotes the Frobenius norm of a matrix, and 
\begin{equation}
	R = [\vc r(t_1)| \vc r(t_2)|\cdots|\vc r(t_M)], \quad \Xi  = [\hat\bxi(t_1)|\hat\bxi(t_2)|\cdots|\hat\bxi(t_M)].
\end{equation}
The solution to~\eqref{eq:Wout_min} is given by
\begin{equation}
	W_{out} = \Xi R^\top \left[RR^\top +\lambda I \right]^{-1}.
\end{equation}

The updates~\eqref{eq:res_state} are designed so that the reservoir state $\vc r(t_i)$ synchronizes with the input time series of the observations $\vc y(t_i)$. As a result, the long-term dynamics of the reservoir should be independent of the initial state $\vc r(t_0)$ of the reservoir. This feature is referred to as the \emph{echo state property}~\cite{Jaeger2001,Maass2002}. To ensure that the echo state property is satisfied, we use a commonly used heuristic, namely that the spectral radius of the reservoir matrix $W_R$ is less than one~\cite{Jaeger2012}. We enforce this property by first defining a sparse matrix $W_R$ whose nonzero entries are randomly drawn from the uniform distribution $U([-0.5,0.5])$. We subsequently rescale this matrix to ensure the desired spectral radius. Similarly, the entries of the input matrix $W_{in}$ and the biases $\vc b$ are drawn at random from the uniform distribution $U([-0.5,0.5])$. The spectral radius of the input weights $W_{in}$ does not play an important role, and therefore no rescaling is needed for $W_{in}$.

To summarize, the observational time series $\vc y(t_i)$ act as an external forcing term that drives the reservoir dynamics~\eqref{eq:res_state}. The reservoir state $\vc r(t_i)$ is mapped to the outputs through the weight matrix $W_{out}$. 
The outputs are then mapped to the kernel vector $\vc z(t_i)$ through the kernel matrix $Z$. Finally, we use this kernel vector in~\eqref{eq:SDEIM}  to arrive at the S-DEIM estimate $\tilde{\vc u}(\vc z(t_i))$. 

Across the examples presented in Section~\ref{sec:num_results}, the training process takes under a second. Once the RC network is trained, the S-DEIM estimate is obtained within milliseconds.
The precise training and inference times are reported in the bottom half of Table~\ref{tab:params}. These numbers are obtained by averaging over ten independent runs on a MacBook Pro with Apple M3 processors. The code is written in Matlab version R2023a. The \emph{reservoir sparsity} refers to the portion of the reservoir matrix $W_R$ that is nonzero.

\section{Numerical Results}\label{sec:num_results}
In this section, we present three numerical examples comparing Q-DEIM and S-DEIM results with fewer sensors than modes ($r<m$). 
The examples are ordered in increasing level of complexity. The first example, the Lorenz-96 system, is a conceptual model of the atmospheric flows governed by a set of nonlinear ODEs. The second example is the Kuramoto--Sivashinsky (KS) equation which is a PDE in one spatial dimension with chaotic dynamics. Finally, we consider the Rayleigh–B\'enard convection (RBC) which is governed by a PDE in two spatial dimensions, modeling the velocity and temperature of a fluid driven by an external temperature gradient.

Table~\ref{tab:params} lists the quantities involved in each example. These include the state space dimension $N$, the number of POD modes $m$, the number of sensors $r$, and the sampling rate $\Delta t$ of the observational data. The number of sensors are taken to be (roughly) half the number of modes. In the PDE examples, the state space dimension $N$ refers to the size of the high-fidelity grid used to discretize the PDE in space.
In order to mimic practical situations where the observational data is polluted with noise, we add a $5\%$ Gaussian noise to every component of the observations $\vc y\in\mathbb R^r$.

For each example, we numerically integrate the model in time from three different random initial conditions. The trajectory corresponding to the first initial condition is used to collect samples for computing the POD modes (cf. Section~\ref{sec:POD}) and the selection matrix (cf. Section~\ref{sec:CPQR}). The second trajectory is used for training the reservoir computing network (cf. Section~\ref{sec:RC_training}). The third trajectory is used for testing purposes. For the training and testing trajectories, the data is sampled every $\Delta t$ time units as reported in Table~\ref{tab:params}.
In the forthcoming sections, only the results corresponding to the testing data are presented. 
Sampling data from three different trajectories demonstrates the robustness and generalizability of S-DEIM. In particular, it shows that the RC network can be trained on one trajectory and applied to a different trajectory as long as the attractor is ergodic.

\begin{table}
	\caption{Parameters for each test case. The bottom half of the table reports the information regarding the RC network.}
	\label{tab:params}
	\centering
	\begin{tabular}{l| c| c| c}
		& \textbf{Lorenz-96} & \textbf{KS} & \textbf{RBC}\\
		 \hhline{=|=|=|=}
		State dimension ($N$) & 40 & 128 & 4257\\
		Number of modes ($m$) & 20 & 15 & 50\\
		Number of sensors ($r$) & 10 & 8 & 25\\
		Sampling rate ($\Delta t$) & 0.25 & 0.2 & 0.5\\ \hline
		Reservoir size ($K$) & 200 & 300 & 500 \\
		Reservoir sparsity  & 10\% & 10\% & 10\% \\
		Training/Inference time (seconds) & 0.2/0.02 & 0.7/0.04 & 0.3/0.02
	\end{tabular}
\end{table}

Finally, in this section, we report the relative error,
\begin{equation}\label{eq:RE}
	\mbox{RE} = \frac{\|\tilde{\vc u}-\vc u\|}{\|\vc u-\overline{\vc u}\|},
\end{equation}
where $\vc u$ is the true state, $\tilde{\vc u}$ is the estimation, and $\overline{\vc u}$ is the empirical mean of the data set from which the POD modes are computed.
Note that in the numerator, the mean $\overline{\vc u}$ cancels out and therefore it is more prudent to subtract the mean from the denominator as well, as in equation~\eqref{eq:RE}. Otherwise, a large mean will result in misleadingly small relative errors.

\subsection{Lorenz-96 model}
The Lorenz-96 equation is a conceptual model of atmospheric flow consisting of $N$ units $\vc u = [u_1,u_2,\cdots,u_N]$ which satisfy the differential equations, 
\begin{equation}
	\dot u_i = -u_i + (u_{i+1}-u_{i-2})u_{i-1}+F,\quad i=1,2,\cdots, N,
	\label{eq:lorenz96}
\end{equation}
with cyclic boundary conditions which imply $u_0 = u_N$, $u_{-1}=u_{N-1}$, and $u_{N+1} = u_1$.
The functions $u_i(t)$ represent some atmospheric quantity, such as potential vorticity, discretized at $N$ points along a circle around the globe at a constant latitude~\cite{Lorenz1996}. The constant $F=4$ represents an external forcing term.
\begin{figure}
	\centering
	\includegraphics[width=\textwidth]{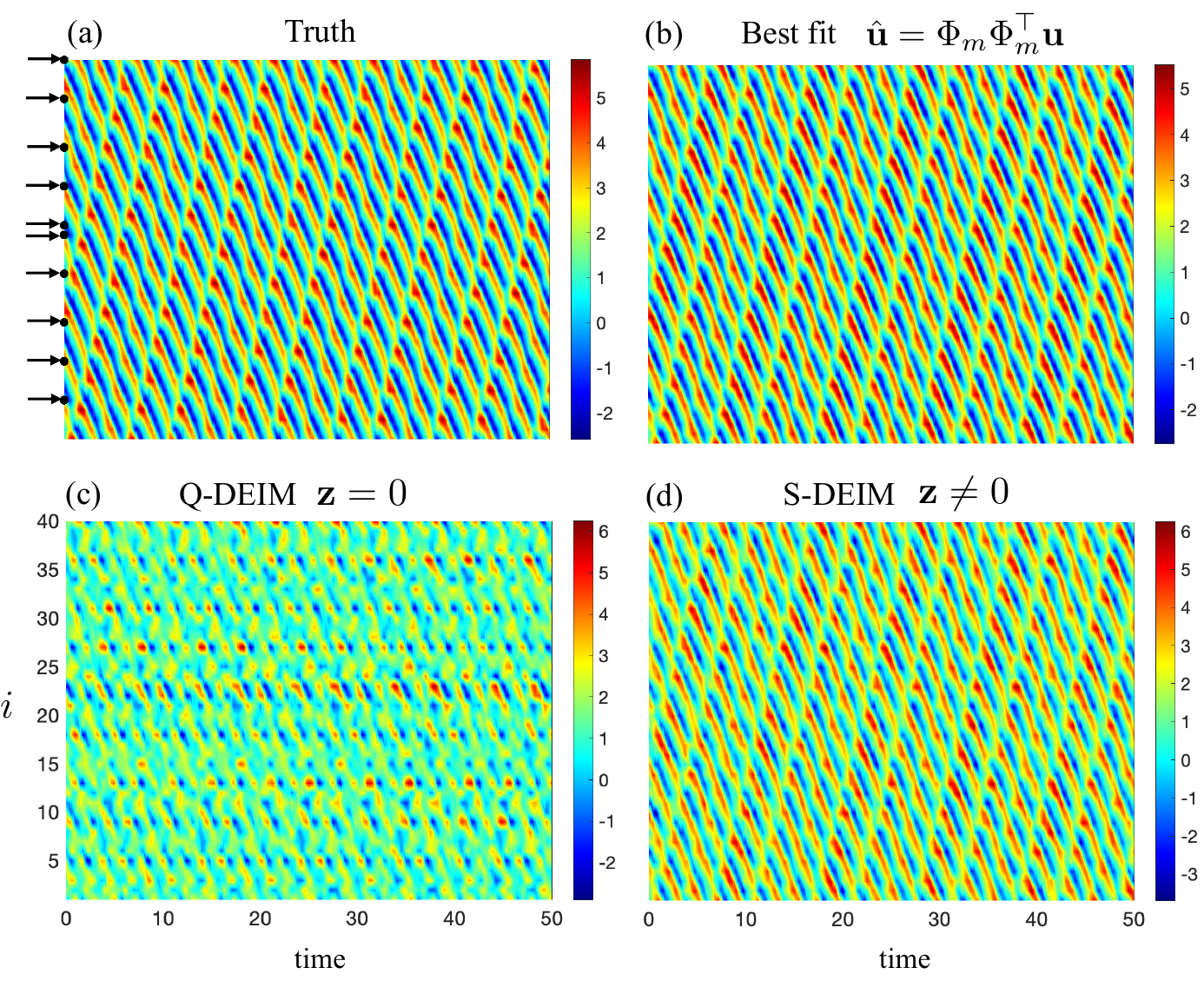}
	\caption{State estimation for Lorenz-96 system. The reconstructions are obtained from $r=10$ sensors and $m=20$ modes}
	\label{fig:lorenz96}
\end{figure}

Figure~\ref{fig:lorenz96}(a) shows a solution of equation~\eqref{eq:lorenz96} with $N=40$ components over the time interval $t\in[0,50]$. This solution exhibits modulated traveling waves. For comparison, figure~\ref{fig:lorenz96}(b) shows the best fit approximation with $m=20$ POD modes, i.e., $\hat{\vc u} = \Phi_m \Phi_m^\top \vc u$. We recall that $\hat{\vc u}$ is the best possible approximation within the truncated POD basis $\Phi_m$. However, this best fit is not computable since the full state $\vc u$ is not known.

In spite of its relatively simple form, the Lorenz-96 system is a challenging problem for POD-based methods. It is well-known that in advection dominated system, a relatively large number of POD modes are required to accurately approximate the state of the system~\cite{Anderson2022}. The traveling waves in the Lorenz-96 system present a similar challenge. As a result, we had to use $m=20$ POD modes (half the size of the system $N=40$) in order for the best fit to reach a relative error of around $20\%$.

For DEIM reconstructions, we use $m=20$ modes and $r=10$ sensors. The sensor locations, marked with arrows in figure~\ref{fig:lorenz96}(a), are determined by the CPQR algorithm, as discussed in Section~\ref{sec:CPQR}.
Figure~\ref{fig:lorenz96}(c) shows the resulting Q-DEIM reconstruction $\tilde{\vc u}(\vc 0)$ which fails to reproduce the main features of the solution, i.e., the modulated traveling waves. 
In contrast, the S-DEIM reconstruction $\tilde {\vc u}(\vc z)$ accurately predicts the traveling wave patterns (panel d), albeit slightly overestimating the amplitude of the oscillations.
In the S-DEIM estimation, the nonzero kernel vector $\vc z$ is the output of a pretrained RC network as outlined in Section~\ref{sec:RC_training}.
This benchmark example highlights two facts: (i) The non-zero kernel vector in S-DEIM plays a crucial role in obtaining more accurate reconstructions and (ii) The pre-trained RNN reliably estimates the optimal kernel vector $\hat{\vc z}(t)$ from the observational time series $\mathcal Y(t)$.

\begin{figure}
	\centering
	\includegraphics[width=\textwidth]{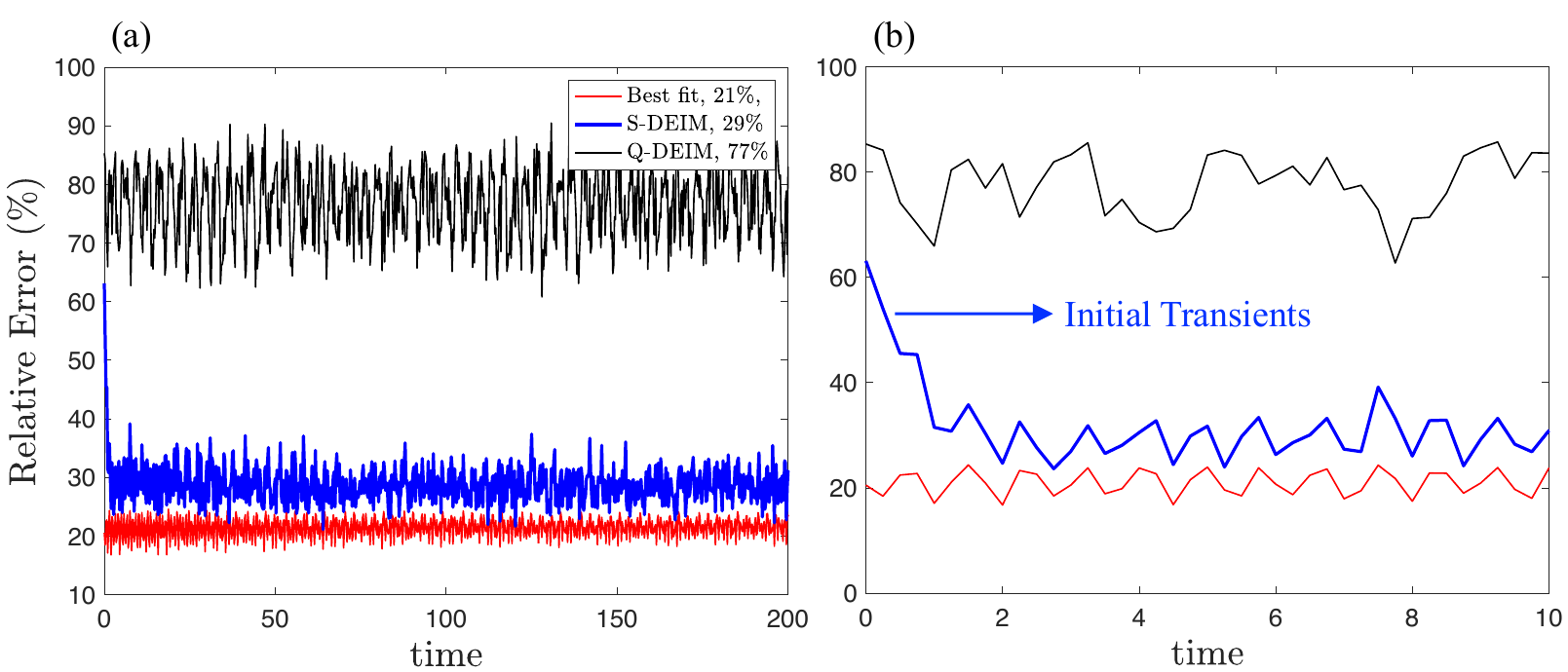}
	\caption{Relative errors for Lorenz-96 system. (a) Relative error as a function of time. The legend reports the mean relative error for each method. (b) A close-up view showing the first 10 time units.}
	\label{fig:lorenz96_err}
\end{figure}

Figure~\ref{fig:lorenz96_err}(a) shows the relative error over the time interval $t\in[0, 200]$. The Q-DEIM error oscillates around 77\%, indicating a poor estimation. In contrast, the S-DEIM error is around 29\% which is close to that of the best possible fit at around 21\%. It is remarkable that although both methods use the same observational data, S-DEIM reduces the reconstruction error by 48\% compared to Q-DEIM.

We recall that the RC network takes some time to synchronize to the observational input data. This transient phase is evident in figure~\ref{fig:lorenz96_err}(b) which shows a close-up view of the relative error. After about two time units, equivalent to $8\Delta t$, the synchronization is complete and the S-DEIM relative error has settled around $29\%$.

\subsection{Kuramoto--Sivashinsky equation}
The Kuramoto--Sivashinsky (KS) system is described by the partial differential equation, 
\begin{equation}
	u_t + uu_x + u_{xx}+u_{xxxx}=0,
\end{equation}
where the solution $u(x,t)$ is defined on the spatial domain $x\in[-L/2,L/2]$ with periodic boundary conditions.
A larger domain size $L$ typically leads to more complicated spatiotemporal dynamics. Here, we set $L=22$ which is known to lead to chaotic dynamics~\cite{Cvitanovic2010}.
We use a standard pseudo-spectral method with $N=128$ collocation points to solve the KS equation~\cite{Trefethen2005}.
\begin{figure}
	\centering
	\includegraphics[width=\textwidth]{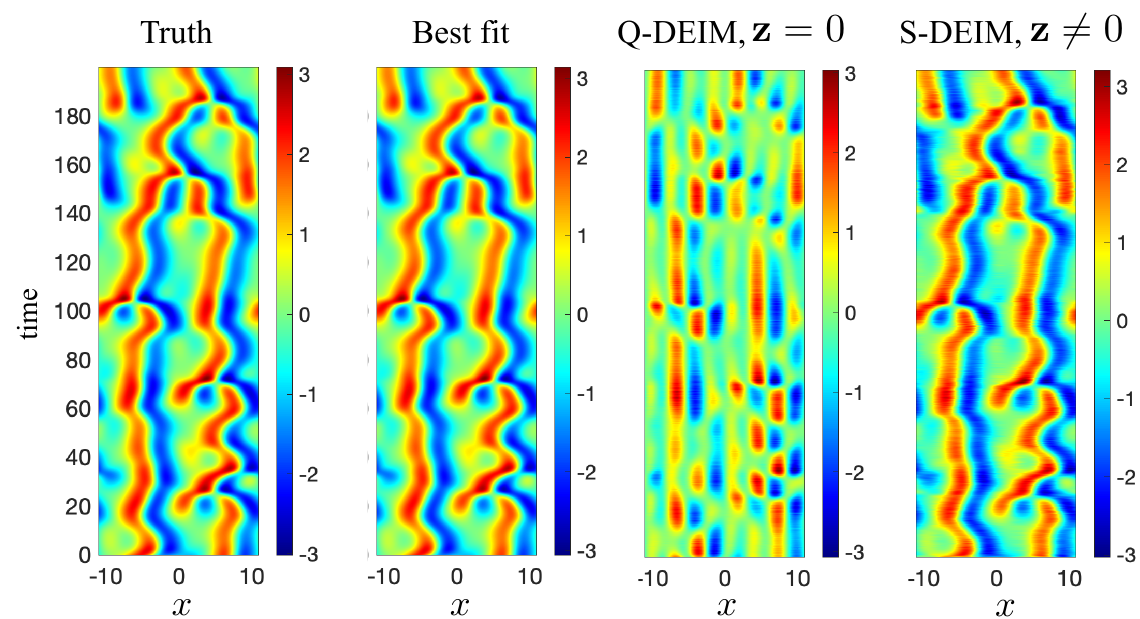}
	\caption{State estimation for the KS system with $m=15$ POD modes and $r=8$ sensors.}
	\label{fig:KS_fields}
\end{figure}

The left panel of figure~\ref{fig:KS_fields} shows a numerical solution up to $t=200$. The solution is saved every $0.2$ time units, resulting in 1001 snapshots of the test data.
For state estimations, we use $m=15$ POD modes and $r=8$ sensors. Figure~\ref{fig:KS_fields} also shows the best fit $\hat{\vc u} = \Phi_m\Phi_m^\top \vc u$, the Q-DEIM estimate $\tilde{\vc u}(\vc 0)$, and the S-DEIM estimate $\tilde{\vc u}(\vc z(t))$. The kernel vector $\vc z(t)$ for S-DEIM is obtained as the output of a pretrained RC network
whose inputs are time series of the observations collected at $r=8$ sensor locations.

\begin{figure}
	\centering
	\includegraphics[width=\textwidth]{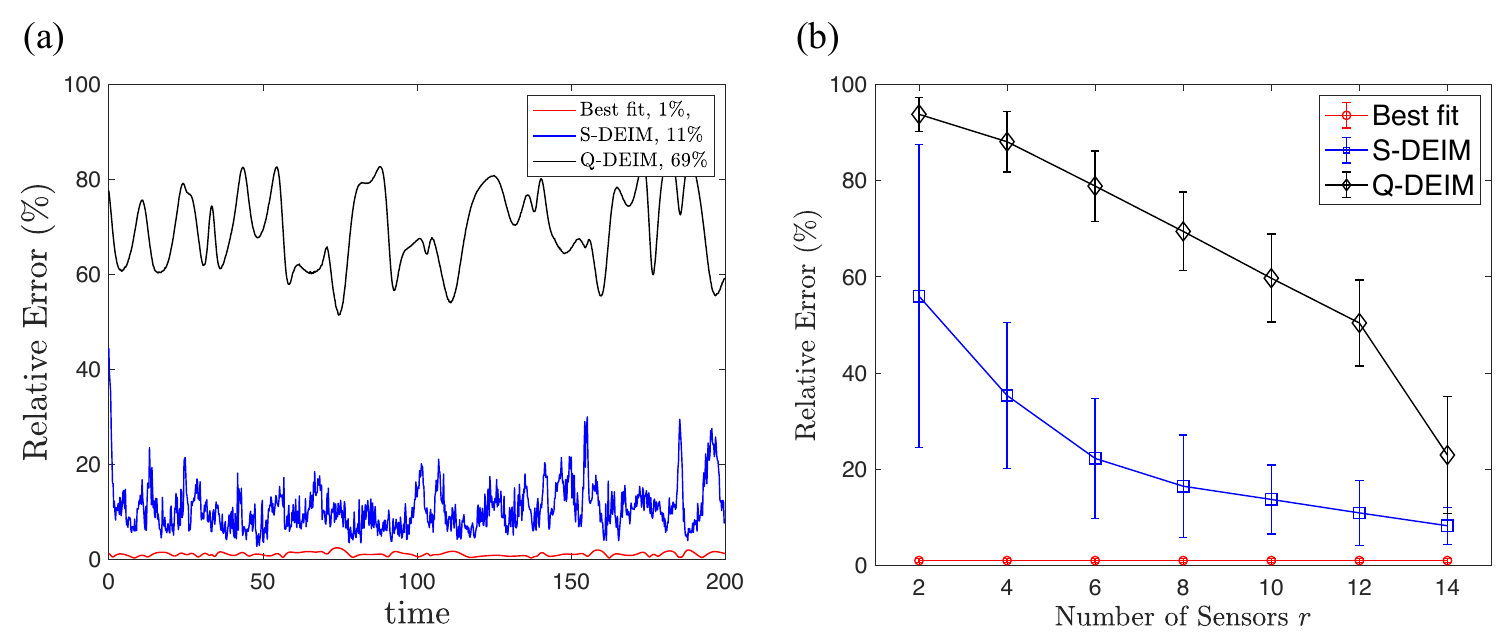}
	\caption{Relative error for the KS equation with $m=15$ POD modes. 
		(a) Relative error as a function of time with $r=8$ sensors. The legend reports the mean relative error for each method.
		(b) Mean relative error as the number of sensors vary. The number of modes is fixed at $m=15$. The error bars mark one standard deviation around the mean.
		Note that the best fit error is independent of the number of sensors.}
	\label{fig:KS_error}
\end{figure}

As shown in figure~\ref{fig:KS_error}(a), the relative error of the best fit is around $1\%$, but recall that the best fit is not computable from the observations. 
The Q-DEIM estimation, which uses $\vc z=\vc 0$, has a mean relative error of about $69\%$. The S-DEIM estimation is much closer to the best fit approximation at a mean relative error of $11\%$.

In figure~\ref{fig:KS_error}(b), we examine the relative error as a function of the number of sensors. We vary the number of sensors between $r=2$ and $r=14$ while keeping the number of POD modes fixed at $m=15$. For each $r$, we compute the mean relative error across 1001 snapshots in the test data. The error bars mark one standard deviation around this mean.
The errors corresponding to both S-DEIM and Q-DEIM decrease as the number of sensors increases. However, the S-DEIM error is always smaller than the Q-DEIM error.  Examining the error bars, it is interesting to observe that the variance of the relative error for S-DEIM decreases as $r$ increases. In contrast, the Q-DEIM variance increases with the number of sensors, indicating an increasing level of uncertainty.

\subsection{Rayleigh–B\'enard convection}
\begin{figure}
	\centering
	\includegraphics[width=\textwidth]{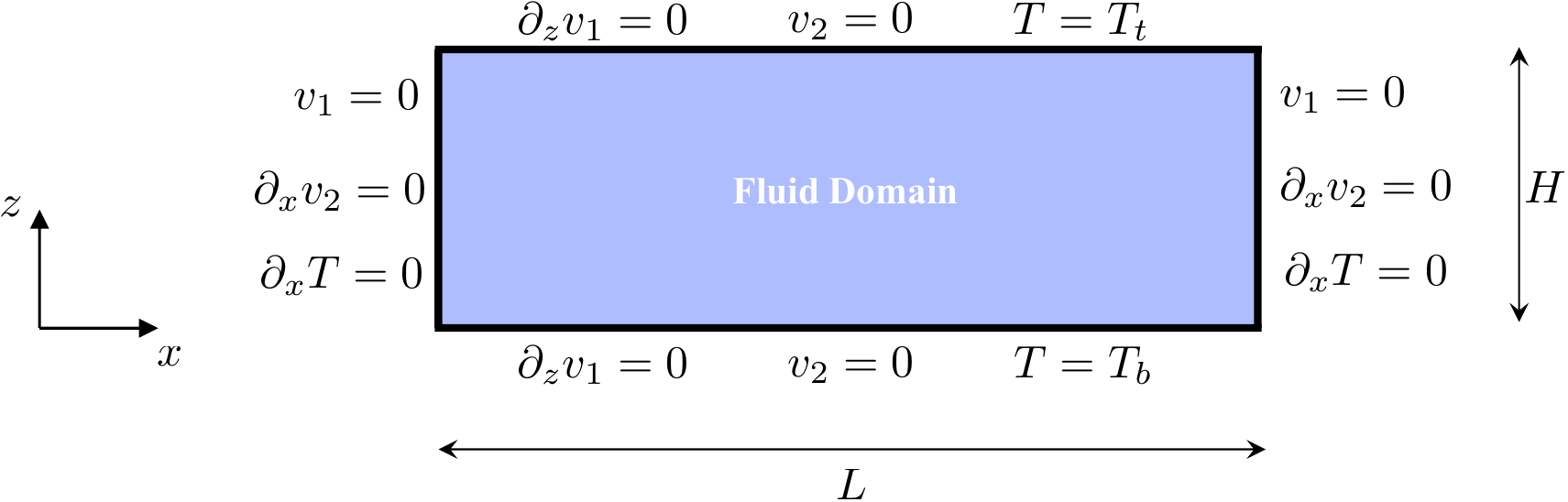}
	\caption{Schematic depiction of the flow geometry for the Rayleigh--B\'enard convection. The boundary conditions are written next to the relevant boundaries.}
	\label{fig:schem_RBC}
\end{figure}

As our last example, we demonstrate the application of S-DEIM on the two-dimensional Rayleigh--B\'enard convection. This problem describes the motion of a fluid inside a rectangular domain $[0,L]\times [0,H]$. As depicted in figure~\ref{fig:schem_RBC}, the bottom wall is kept at a relatively high temperature $T_b$ and the top wall is maintained at a lower temperature $T_t$, where $T_b>T_t$. The fluid is heated by the bottom wall and rises up towards the top wall. The fluid is subsequently cooled down and descends towards the bottom wall again. This process leads to convection cells which can exhibit chaotic dynamics~\cite{Schumacher2018,Kooloth2021}.

In dimensionless variables, the governing equations for this system are given by 
\begin{subequations}
\begin{equation}
\partial_t \vc v + \vc v\cdot\nabla \vc v = -\nabla p + T\vc e_z +\sqrt{\frac{Pr}{Ra}}\Delta \vc v,\quad \nabla\cdot\vc v = 0,
\end{equation}
\begin{equation}
	\partial_t T + \vc v\cdot\nabla T = \frac{1}{\sqrt{PrRa}}\Delta T,
\end{equation}
\end{subequations}
where $T(x,z,t)$ denotes the fluid temperature, $p(x,z,t)$ denotes pressure, and $\vc v(x,z,t)$ is the fluid velocity field. The unit vector $\vc e_z$ points in the opposite direction of gravity. This model has two dimensionless parameters: Rayleigh number $Ra$ and Prandtl number $Pr$. Throughout this section, we set $Ra=10^5$ and $Pr=10$. For the domain size, we use $L=4$ and $H=1$.
The relevant boundary conditions are listed in figure~\ref{fig:schem_RBC}. In our numerical simulations, we set $T_b=1$ and $T_t=0$ for the wall temperatures. A Fourier pseudo-spectral method is used for discretizing the PDEs in space and a 4th order Runge-Kutta scheme is used for temporal discretization~\cite{Verma2013,Hilliard2024b}.
\begin{figure}
	\centering
	\includegraphics[width=.5\textwidth]{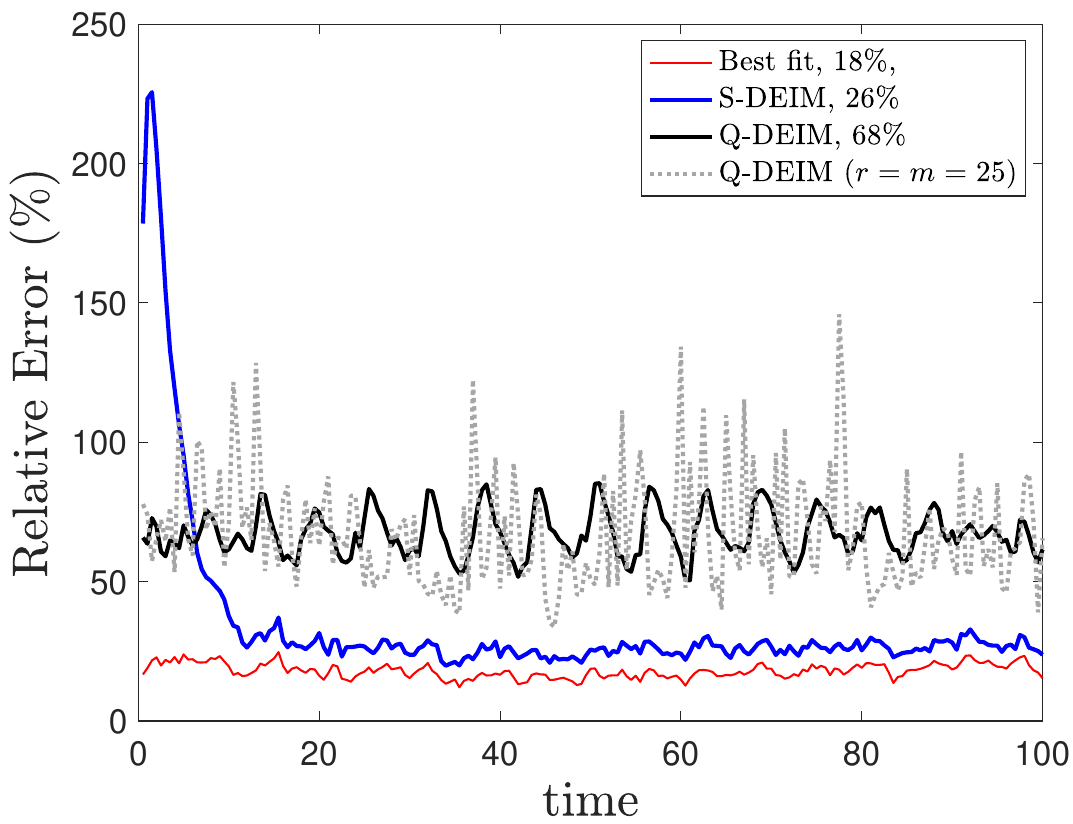}
	\caption{Relative error for the Rayleigh--B\'enard convection using S-DEIM and Q-DEIM with $r=25$ sensors and $m=50$ modes. The best fit (not computable from observations) is also shown for reference. The numbers in the legend report the mean relative error. The dotted gray curve corresponds to the Q-DEIM results with $r=m=25$.}
	\label{fig:RBC_error}
\end{figure}

We seek to estimate the temperature field $T(x,z,t)$ from its sparse observations. The velocity field $\vc v(x,z,t)$ is not observed or used in these estimations.
Figure~\ref{fig:RBC_error} shows the relative errors of state estimations using $r=25$ sensors and $m=50$ POD modes. The Q-DEIM estimates, which use the trivial kernel vector $\vc z=\vc 0$, have a mean relative error of approximately $68\%$. The S-DEIM reconstruction, which uses an RNN to estimate the optimal kernel vector, achieves much smaller errors. Recall that RNNs require some time to synchronize with the input data. As such, the S-DEIM errors are initially large. But within approximately $15$ time units, equivalent to $30\Delta t$, the error decays towards its mean relative error of approximately $26\%$.

For reference, in figure~\ref{fig:RBC_error}, we also show the relative error corresponding to the best approximation $\hat{\vc u}=\Phi_m\Phi_m^\top\vc u$. Recall that this is the best possible approximation within the truncated POD basis $\Phi_m$, however it is not computable from observational data. The mean relative error of the best approximation is about $18\%$, showing that the S-DEIM approximation is reasonably close to the best possible approximation. 

\begin{figure}
	\centering
	\includegraphics[width=\textwidth]{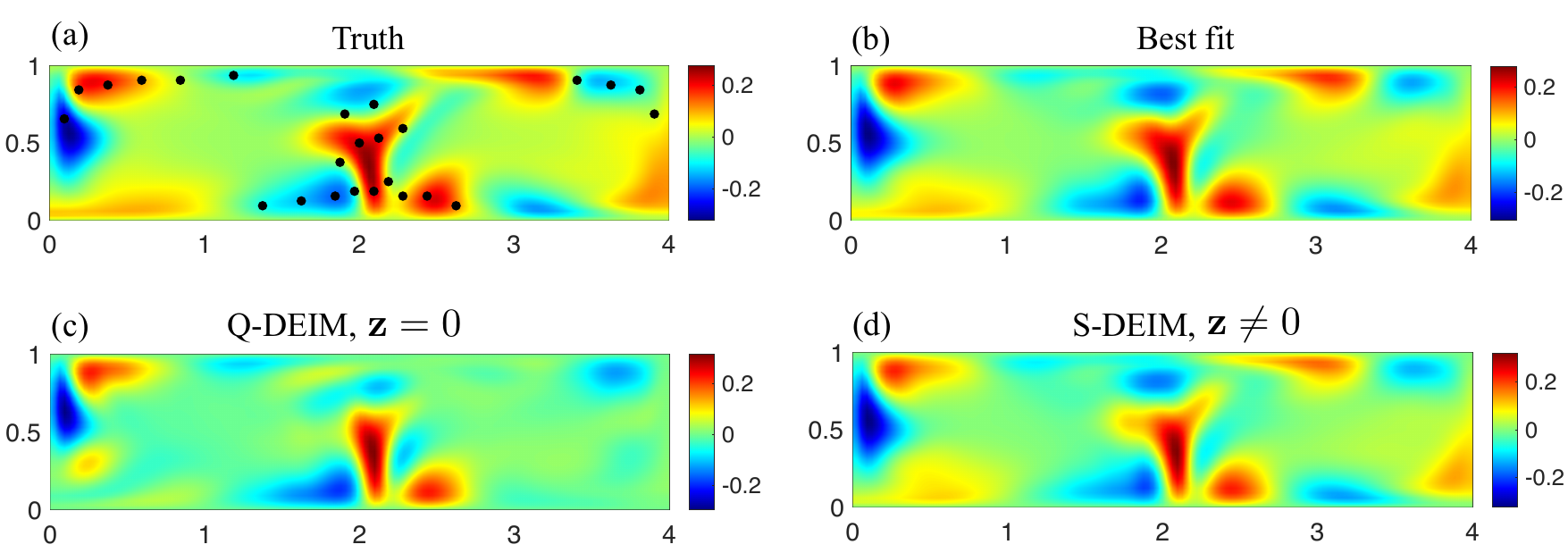}
	\caption{State estimation for the Rayleigh--B\'enard convection with $r=25$ sensors and $m=50$ modes. Each panel shows the mean-removed temperature, $T(x,z,t)-\overline{T}(x,z)$, at time $t=100$. The black circles in panel (a) mark the sensor locations.}
	\label{fig:RBC_T}
\end{figure}

This is also demonstrated in figure~\ref{fig:RBC_T}, which shows the state estimations at time $t=100$. In order to see the differences more clearly, we show the centered data $T(x,z,t)-\overline{T}(x,z)$ where the time-average $\overline T(x,z)$ is subtracted from the temperature.
Figure~\ref{fig:RBC_T}(a) shows the true state to be estimated. This panel also shows the location of $25$ sensors obtained from the CPQR algorithm (black circles).

Figure~\ref{fig:RBC_T}(b) shows the best fit approximation. The S-DEIM approximation (panel d) closely resembles the best fit and captures most features of the true state. In contrast, the Q-DEIM estimate (panel c) fails to reproduce some finer features of the state.

We conclude this section by commenting on the number of modes $m$ used in the Q-DEIM reconstruction. The original version of Q-DEIM was developed under the assumption that the number of sensors and modes coincide, $m=r$~\cite{Sorensen2010,Drmac2016}. One can argue that if a limited number of sensors $r$ is available, the number of modes should be reduced to match $r$.
However, reducing the number of modes to $m=r$ does not necessarily improve the Q-DEIM estimate. For instance, the dotted gray line in figure~\ref{fig:RBC_error} shows the Q-DEIM estimation error with $m=r=25$. We observe that reducing the number of modes to match the available number of sensors does not improve the estimation of the temperature field in the Rayleigh--B\'enard convection. Therefore, a larger number of POD modes are needed to approximate the state of this system, underscoring the need for tackling the underdetermined case where $m>r$, as S-DEIM seeks to do.

\subsection{Comparison with direct learning}
Recall the best fit approximation~\eqref{eq:bestFit} which corresponds to the expansion coefficients $\hat\balpha(t) = \Phi_m^\top \vc u(t)$.
Although these coefficients are not directly computable from observations, we can nonetheless try to estimate them from observational time series $\mathcal Y(t)$ using an RNN.
More precisely, we can train an RNN whose inputs are the observational time series $\mathcal Y(t)$ (similar to S-DEIM) and its outputs are the best fit coefficients $\hat\balpha(t)$ (unlike S-DEIM where the outputs are the optimal kernel vector $\hat{\bxi}(t)$). We refer to this approach as the \emph{direct learning method} and in this section compare its performance to S-DEIM.

The direct method has two drawbacks when compared to S-DEIM. First, the dimension of the RNN output space is larger for the direct method. For the direct method, the outputs $\hat\balpha$ are $m$-dimensional, whereas in S-DEIM the outputs $\hat{\bxi}$ are $(m-r)$-dimensional. This may render the training of the RNN for the direct method slightly more challenging. But, the more important issue is the loss of the interpolation property~\eqref{eq:interp}. Owing to the interpolation property, S-DEIM estimates agree exactly with the observations. In contrast, there is no such guarantee for the direct method. 

\begin{figure}
	\centering
	\includegraphics[width=\textwidth]{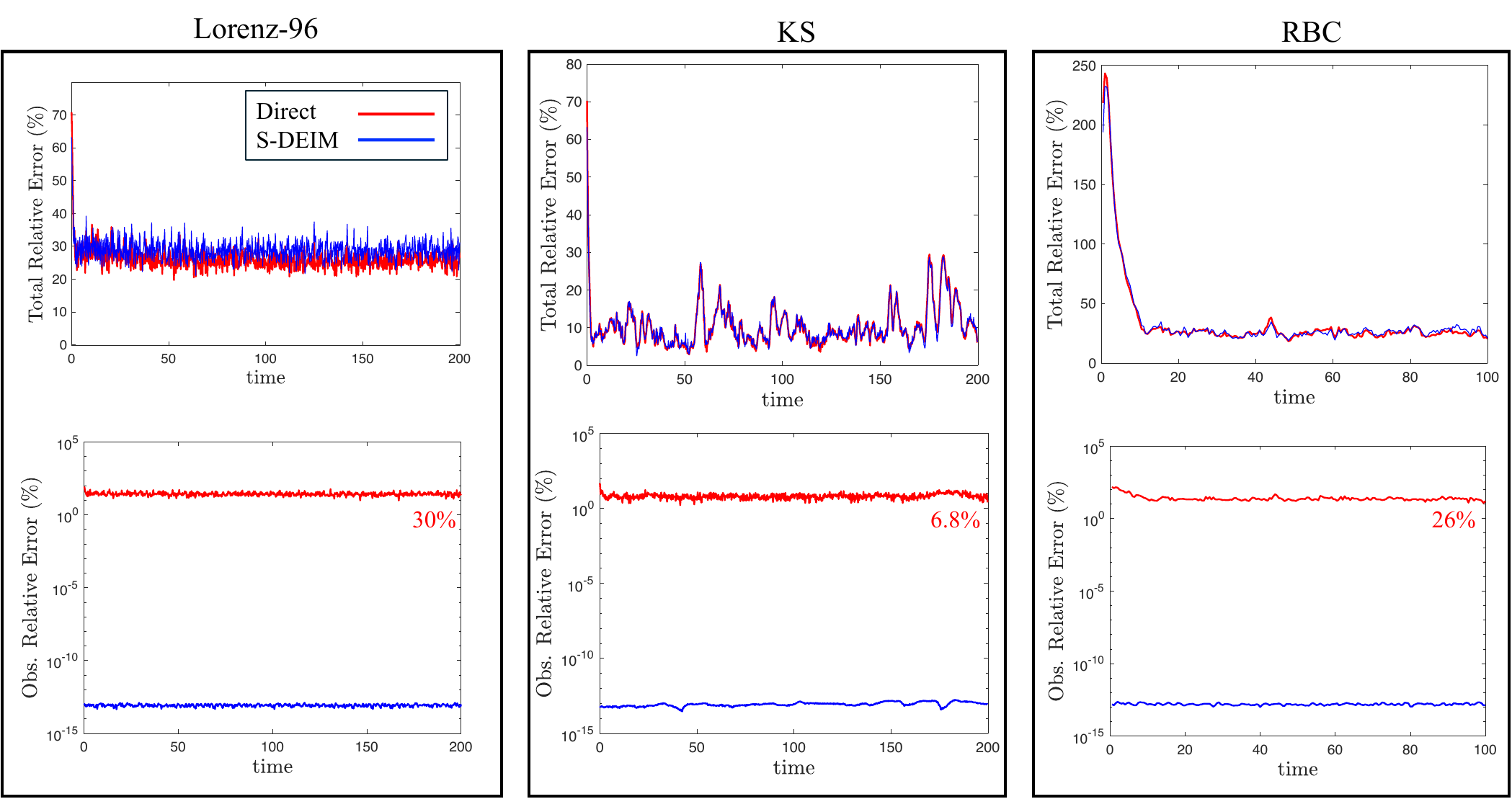}
	\caption{Comparison of S-DEIM and direct learning. The top row shows the total relative error~\eqref{eq:RE}. The bottom row shows the relative observation error. The numbers report the time-average of the observation errors for the direct method.}
	\label{fig:SDEIMvsDirect_comp}
\end{figure}

In spite of these drawbacks, the direct method should be theoretically more accurate since it seeks to estimate the best fit. However, as we show in figure~\ref{fig:SDEIMvsDirect_comp} this higher accuracy does not materialize in practice. The top row of figure~\ref{fig:SDEIMvsDirect_comp} shows the total relative error~\eqref{eq:RE} for the direct learning method and S-DEIM. The total relative errors for both methods are remarkably similar, across all three examples: Lorenz-96, KS, and RBC systems. 
Here, to  perform a fair comparison, we used the same RC architecture for the direct learning method and for S-DEIM, except for the output dimensions. Furthermore, both methods are trained on the same input data and also tested on identical datasets. Figure~\ref{fig:SDEIMvsDirect_comp} shows the results on the test data.

As mentioned earlier, S-DEIM is guaranteed to satisfy the interpolation property. This is shown in the bottom row of figure~\ref{fig:SDEIMvsDirect_comp} where the relative observation error $\| S_r^\top \tilde{\vc u}-\vc y\|/\|\vc y\|$ is plotted. The S-DEIM observation error is about $10^{-13}\%$ for all three systems; it is not exactly zero because of numerical round-off errors. 
The observation error for the direct method is at least thirteen orders of magnitude larger, ranging between 6.8-30\% on average.

In summary, we observe that the direct method not only loses the interpolation property but also fails to improve the estimation error appreciably. As a result, we do not recommend using the direct learning method instead of S-DEIM.

\section{Conclusions}\label{sec:conc}
We proposed an model-free method for estimating the state of a spatiotemporal dynamical system from incomplete observations. Our method relies on Sparse Discrete Empirical Interpolation Method (S-DEIM) which comprises two terms. The first term can immediately be evaluated using the observations. The second term involve the kernel vector of a linear operator which cannot be estimated from instantaneous observational data. Here, we developed a machine learning method to estimate this optimal kernel vector. More precisely, we train a recurrent neural network whose inputs are the time series of the observations (including past history) and its output estimates the optimal kernel vector of S-DEIM.

We demonstrated on three numerical examples that S-DEIM, together with the RNN training, returns state estimations with nearly optimal kernel vectors. The resulting S-DEIM estimations are 42-58\% more accurate than the Q-DEIM estimates. Furthermore, the computational cost of the proposed method is exceedingly low. In each numerical example, the training process took only a few seconds. Using the trained networks to obtain a state estimation took only a fraction of a second.

With certain adjustments, our method is scalable to large scale systems, such as three-dimensional turbulent flows, where the state space dimension $N$ is extremely large. The bulk of the computational cost is associated with the offline computations, as depicted in figure~\ref{fig:schem_pipeline}. Within this stage, producing the high-fidelity experimental or simulation data and the subsequent POD computation are the most resource intensive steps. For instance, in large scale applications, it is computationally infeasible to compute POD modes from a direct singular value decomposition. Instead, one should use low cost and low memory methods such as randomized matrix decomposition~\cite{Tropp2011} or tensor-based methods~\cite{Farazmand2023}. It is worth mentioning that the training phase is not significantly impacted by the size $N$ of the state space. Even in large scale systems, the input and output dimensions of the RNN are determined by the number of sensors $r$ and the number of modes $m$, which are relatively small.

Here, we used a simple RNN architecture, namely the reservoir computing network. Future work will study different types of RNNs such as LSTM networks. A more fundamental open problem concerns sensor placement for S-DEIM. Here, we used the CPQR decomposition for sensor placement which is tailored towards Q-DEIM where the trivial kernel vector $\vc z=\vc 0$ is used. For S-DEIM, which uses a nonzero kernel vector, optimal sensor locations will likely differ. As such, new sensor placement methods, tailor-made for S-DEIM, are needed. 

\subsubsection*{Data and code} The codes for this work can be found at \href{https://github.com/mfarazmand/S-DEIM.git}{https://github.com/mfarazmand/S-DEIM.git} and the corresponding data is available at\\
\href{https://doi.org/10.5281/zenodo.15832022}{https://doi.org/10.5281/zenodo.15832022}.

\subsubsection*{Acknowledgments}
This work was supported by the National Science Foundation, the Algorithms for Threat Detection (ATD) program, through the award DMS-2220548.
I am grateful to Zack Hilliard for help with setting up the numerical simulations for the Rayleigh--B\'enard convection. I would also like to thank Prof. Arvind Saibaba for his comments on an earlier version of this paper.


\end{document}